\documentclass[11pt]{amsart}
\usepackage[english]{babel}
\usepackage[latin1]{inputenc}
\usepackage{amstext}
\usepackage{amsmath}
\usepackage{amsfonts}
\usepackage{amssymb}
\usepackage{amsthm}
\usepackage[all]{xy}
\xyoption{arc}
\CompileMatrices
\usepackage{graphicx}

\newtheorem{thm}{Theorem}[section]
\newtheorem{cor}[thm]{Corollary}
\newtheorem{lem}[thm]{Lemma}
\newtheorem{prop}[thm]{Proposition}
\newtheorem{dfn}[thm]{Definition}
\newtheorem{rem}[thm]{Remark}

\newenvironment{pr}[1][Proof]{\noindent\textbf{#1.} }{\ \rule{0.5em}{0.5em}}

\newcommand{\SO}{{\mathcal{O}}}
\newcommand{\PP}{\mathbb{P}}

\newcommand{\CC}{\mathbb{C}}

\newcommand{\Hom}{\operatorname{Hom}}

\newcommand{\surj}{\twoheadrightarrow}

\newcommand{\too}{\longrightarrow}
\newcommand{\rk}{\operatorname{rk}}

\newcommand{\GL}{\operatorname{GL}}
\newcommand{\SL}{\operatorname{SL}}

\selectfont

\numberwithin{equation}{section}

\begin{document}

\title{A GIT interpretation of the Harder-Narasimhan filtration}
\author[T. G\'omez, I. Sols, A. Zamora]
{Tom\'as L. G\'omez, Ignacio Sols and Alfonso Zamora}


\address{TG, AZ: Instituto de Ciencias Matem\'aticas (CSIC-UAM-UC3M-UCM),
Nicol\'as Cabrera 13-15, Campus Cantoblanco UAM, 28049 Madrid,
Spain}

\address{IS, AZ: Departamento de \'{A}lgebra, Facultad de Matem\'aticas, Universidad Complutense de
Madrid, 28040 Madrid, Spain}

\email{tomas.gomez@icmat.es, isols@mat.ucm.es, alfonsozamora@icmat.es}

\begin{abstract}
An unstable torsion free sheaf on a smooth projective variety gives a
GIT unstable point in certain Quot scheme.  To a GIT unstable point,
Kempf associates a ``maximally destabilizing'' 1-parameter subgroup,
and this induces a filtration of the torsion free sheaf. We show that
this filtration coincides with the Harder-Narasimhan filtration.
\end{abstract}

\maketitle

\section*{Introduction}

Let $X$ be a smooth complex projective variety, and let $\SO_X(1)$ be an ample
line bundle on $X$. If $E$ is a coherent sheaf on $X$, let $P_E$ be its Hilbert
polynomial with respect to $\SO_X(1)$, i.e.,
$P_E(m)=\chi(E\otimes \SO_X(m))$. If $P$ and $Q$ are polynomials, we
write $P\leq Q$ if $P(m)\leq Q(m)$ for $m\gg 0$.

A torsion free sheaf $E$ on
$X$ is called semistable if for all proper subsheaves $0\neq F\subset E$,
$$
\frac{P_F}{\rk F} \leq \frac{P_E}{\rk E} \; .
$$
If it is not semistable, it is called unstable, and it has a canonical
filtration:

Given a torsion free sheaf $E$, there exists a unique filtration
$$0 =E_0 \subset E_{1} \subset E_{2} \subset \cdots \subset E_{t} \subset
E_{t+1}=E\; ,$$ which satisfies the following properties, where
$E^{i}:=E_{i}/E_{i-1}$:
 \begin{enumerate}
   \item Every $E^{i}$ is semistable
   \item The Hilbert polynomials verify
   $$\frac{P_{E^{1}}}{\rk E^{1}}>\frac{P_{E^{2}}}{\rk E^{2}}>\ldots>\frac{P_{E^{t+1}}}{\rk E^{t+1}}$$
 \end{enumerate}
This filtration is called the \emph{Harder-Narasimhan filtration} of
$E$ (\cite[Theorem 1.3.6]{HL2}).

We will briefly describe the construction of the moduli space for these objects.
This is originally due to Gieseker for surfaces, and it was
generalized to higher dimension by Maruyama (\cite{Gi,Ma}).
To construct the moduli space of torsion free sheaves with fixed
Hilbert polynomial $P$, we choose a suitably large integer $m$ and
consider the Quot scheme parametrizing quotients
\begin{equation}
\label{quot}
V\otimes \SO_X(-m) \too E
\end{equation}
where $V$ is a fixed vector space of dimension $P(m)$ and $E$ is a
sheaf with $P_E=P$.
The Quot scheme has a canonical action by
$\SL(V)$. Gieseker (c.f. \cite{Gi}) gives a linearization of this action
on a certain ample line bundle, in order to use Geometric
Invariant Theory to take the quotient by the action.
The moduli space of semistable sheaves is obtained
as the GIT quotient.

Let $E$ be a torsion free sheaf which is unstable. Choosing $m$
large enough (depending on $E$),
and choosing an isomorphism $V\cong H^0(E(m))$, we
obtain a quotient as in \eqref{quot}. The corresponding point in
the Quot scheme will be GIT unstable. By the Hilbert-Mumford
criterion, there will be a 1-parameter subgroup of $\SL(V)$ which
``destabilizes'' the point. Among all these 1-parameter subgroups,
Kempf (c.f. \cite{Ke}) shows that there is a conjugacy class of ``maximally
destabilizing'' 1-parameter subgroups, all of them giving a unique
weighted filtration of $V$. This filtration induces a sheaf
filtration of $E$. In principle, this filtration will depend on
the integer $m$ but we show that it stabilizes for $m\gg 0$, and
we call it the \emph{Kempf filtration of $E$}. In this article, we
show that the Kempf filtration of an unstable torsion free sheaf
$E$ coincides with the Harder-Narasimhan filtration.

If $X$ is a curve, this result is in the Master's Thesis 
of Alfonso Zamora in 2009 (c.f. \cite{Za2}).
The use of Kempf's ideas is
already there, and in this article the method is extended and refined to 
obtain the result in higher dimension.

The equality between the Harder-Narasimhan filtration and the
Kempf filtration for torsion free sheaves has independently
been studied by
Hoskins and Kirwan (c.f. \cite{HK}) in the stratification of the
Quot in Harder-Narasimhan types. The difference with our approach
is that they use the existence of the Harder-Narasimhan
filtration, having fixed the Harder-Narasimhan type for each
stratum, whereas we prove that the Kempf filtration is independent of
$m$ if $m$ is large enough without using the Harder-Narasimhan
filtration. In other words, our method gives a different proof of
the existence of the Harder-Narasimhan filtration, and in principle our
method could be used to define the Harder-Narasimhan filtration
(using the Kempf filtration for $m$ large) 
in a moduli problem where there is
still no Harder-Narasimhan filtration known. This is the case of
\cite{Za1}, where a similar construction is developed for rank $2$
tensors.

In fact, the motivation for this work was to give a general
procedure to obtain the Harder-Narasimhan filtration for any 
moduli problem constructed with GIT, using the Kempf filtration.
The main difficulty we found was to prove that the filtration that
we obtain is independent of the integer $m$, once this is large enough.

One referee suggested to use the fact,
proved in \cite{RR}, that the limiting point of the one parameter
subgroup $\lambda$ given by Kempf is semistable with respect to the
induced action of the reductive centralizer of $\lambda$. In the case
of torsion free sheaves this will imply that, for suitably large $m$, the
successive quotients of the filtration induced on $E$ by $\lambda$ are
semistable, and this is one of the properties of the Harder-Narasimhan
filtration. It would be interesting to study if this approach can be
used in other moduli problems.

We hope that our approach will be useful to find a Harder-Narasimhan
filtration in situations where it is still not defined.


If we replace Hilbert polynomials with degrees, the notion of
semistability becomes $\mu$-semistability
(also known as slope semistability) and we obtain the
$\mu$-Harder-Narasimhan filtration. In \cite{Br,BT},
Bruasse and Teleman give a gauge-theoretic interpretation of
the $\mu$-Harder-Narasimhan filtration for torsion free sheaves
and for holomorphic pairs over holomorphic curves, where stability and $\mu$-stability do coincide. 
They also use Kempf's ideas, but generalizing them to the setting of the gauge group, to show analogous correspondences in the complex geometry framework.

{\bf Acknowledgments.} We thank Francisco Presas for discussions.
This work was funded by the grant MTM2010-17389 and ICMAT Severo Ochoa
project SEV-2011-0087 of the Spanish
Ministerio de Econom{\'\i}a y Competitividad. 
A. Zamora was supported by a
FPU grant from the Spanish Ministerio de Educaci\'on. Finally A.
Zamora would like to thank the Department of Mathematics at
Columbia University, where part of this work was done, for
hospitality. This work is part of A. Zamora's Ph.D. thesis(c.f.
\cite{Za3}).

\section{A theorem by Kempf}

Following the usual convention, whenever
``(semi)stable'' and ``$(\leq)$'' appear in a sentence, two
statements should be read: one  with
``semistable'' and ``$\leq$'' and another
with ``stable'' and ``$<$''.

Let $X$ be a  smooth complex projective variety of dimension $n$
endowed with a fixed polarization $\SO_X(1)$. A torsion free sheaf
$E$ on $X$ is said to be (semi)stable if for all non zero proper subsheaves
$F$
\begin{equation}
  \label{defss}
\frac{P_F}{\rk F} \; (\leq) \; \frac{P_E}{\rk E} \; .
\end{equation}
We will recall Gieseker's
construction  (c.f. \cite{Gi}) of the moduli space of semistable torsion free sheaves
with fixed Hilbert polynomial $P$ and fixed determinant $\det(E)\cong
\Delta$\ .

A coherent sheaf is called $m$-regular if $h^i(E(m-i))=0$ for all
$i>0$.
\begin{lem}
\label{mregularity}
If $E$ is $m$-regular then the following holds
\begin{enumerate}
\item $E$ is $m'$-regular for $m'>m$
\item $E(m)$ is globally generated
\item For all $m'\geq 0$ the following homomorphisms are surjective
$$
H^0(E(m))\otimes H^0(\SO_X(m'))\too H^0(E(m+m')) \; .
$$
\end{enumerate}
\end{lem}

Let $m$ be a suitable large integer, so that $E$ is $m$-regular
for all semistable $E$ (c.f. \cite[Corollary 3.3.1 and Proposition
3.6]{Ma}). Let $V$ be a vector space of dimension $p:=P(m)$. Given an
isomorphism $V\cong H^0(E(m))$ we obtain a quotient
$$
q:V\otimes \SO_X(-m) \surj E\; ,
$$
hence a homomorphism
$$
Q:\wedge^r V \cong \wedge^r H^0(E(m)) \too H^0(\wedge^r(E(m)))
\cong H^0(\Delta(rm))=:A\; ,
$$
and points
$$
Q\in \Hom(\wedge^r V , A) \qquad\overline{Q} \in
\PP(\Hom(\wedge^r V , A) )\; ,
$$
where $Q$ is well defined up to a
scalar because the isomorphism $\det(E)\cong \Delta$
is well defined up to a scalar, and hence
$\overline{Q}$ is a well defined point.
The point $\overline{Q}$ depends on $E$ and on the chosen isomorphism
$V\cong H^0(E(m))$. To get rid of the choice of isomorphism, we have
to take the quotient by the canonical action of $\GL(V)$. Since an
homothecy acts trivially, we might as well take the quotient by
$\SL(V)$.

A weighted filtration $(V_\bullet,n_\bullet)$
of $V$ is a filtration
\begin{equation}
\label{filtV} 0\subset V_1 \subset V_2 \subset \;\cdots\; \subset V_t \subset V_{t+1}=V,
\end{equation}
and rational numbers $n_1,\, n_2,\ldots , \,n_t > 0$.
To a weighted filtration we associate a vector of $\mathbb{Q}^p$
defined as $\Gamma=\sum_{i=1}^{t}n_i \Gamma^{(\dim V_i)}$
where
\begin{equation}
\label{semistandard}
\Gamma^{(k)}:=\big( \overbrace{k-p,\ldots,k-p}^k,
 \overbrace{k,\ldots,k}^{p-k} \big)
\qquad (1\leq k < p)
\, .
\end{equation}
Hence, the vector is of the form
$$\Gamma=(\overbrace{\Gamma_1,\ldots,\Gamma_1}^{\dim V^1},
\overbrace{\Gamma_2,\ldots,\Gamma_2}^{\dim V^2},
\ldots,
\overbrace{\Gamma_{t+1},\ldots,\Gamma_{t+1}}^{\dim V^{t+1}}) \; ,$$
where $V^i=V_i/V_{i-1}$. Giving the numbers $n_1,\ldots,n_t$ is clearly equivalent
to giving the numbers $\Gamma_1,\ldots,\Gamma_{t+1}$
because
$$n_i=\frac{\Gamma_{i+1}-\Gamma_i}{p}
\qquad \text{and} \quad \sum_{i=1}^{t+1}\Gamma_i\dim V^i = 0$$
A 1-parameter subgroup
of $\SL(V)$ (which we denote in the following by 1-PS)
is a non-trivial homomorphism
$\CC^*\to \SL(V)$. To a 1-PS we associate a weighted
filtration as follows. There is a basis $\{e_1,\ldots,e_p\}$ of
$V$ where it has a diagonal form
$$
t\mapsto \operatorname{diag} \big(
t^{\Gamma_1},\ldots,t^{\Gamma_1},
t^{\Gamma_2},\ldots,t^{\Gamma_2}, \ldots ,
t^{\Gamma_{t+1}},\ldots,t^{\Gamma_{t+1}} \big)
$$
with
$\Gamma_1<\cdots<\Gamma_{t+1}$. Let
$$
0\subset V_1 \subset \cdots \subset V_{t+1}=V
$$
be filtration obtained by calling $V_{1}\subset V$ the vector subspace generated by those vectors of the basis of $V$
associated to exponents $\Gamma_{1}$, $V_{2}\subset V$ generated for those associated to exponents $\Gamma_{1}$ and $\Gamma_{2}$, and so on. 
Note that two 1-PS give the same
filtration if and only if they are
conjugate by an element of the parabolic subgroup of
$\SL(V)$ defined by the filtration.

The basis $\{e_1,\ldots,e_p\}$, together with a
basis $\{w_j\}$ of $A$, induces a basis of
$\Hom(\wedge^r V,A)$ indexed in a natural way by tuples
$(i_1,\ldots,i_r,j)$ with $i_1<\cdots <i_{r}$, and the
coordinate corresponding to such an index is acted by
the 1-PS as:
$$
Q_{i_1,\cdots,i_r,j} \mapsto t^{\Gamma_{i_1}+\cdots+\Gamma_{i_r}}
Q_{i_1,\cdots,i_r,j}
$$
The coordinate $(i_1,\ldots,i_r,j)$ of the point corresponding to
$E$ is non-zero if and only if the evaluations of the sections
$e_1,\ldots,e_r$ are linearly independent for generic $x\in X$.
Therefore, the ``minimal relevant weight'' which has to be
calculated to apply Hilbert-Mumford criterion for GIT stability is
\begin{eqnarray}
  \label{eq:muleft}
\mu(\overline{Q},V_\bullet,n_\bullet)&=&\min \{\Gamma_{i_1}+\cdots+\Gamma_{i_r}: \,
Q_{i_1,\ldots ,i_r,j}\neq 0 \} \notag\\
&=&\min \{\Gamma_{i_1}+\cdots+\Gamma_{i_r}: \,
Q(e_{i_1}\wedge\cdots \wedge e_{i_r})\neq 0 \} \notag\\
&=&\min \{\Gamma_{i_1}+\cdots+\Gamma_{i_r}: \,
e_{i_1}(x),\ldots,e_{i_r}(x)  \\
& & \qquad\text{linearly independent for generic
  $x\in X$}\}   \notag
\end{eqnarray}

After a short calculation (originally due to Gieseker) we obtain
\begin{equation}
\label{mrw}
\mu(\overline{Q},V_\bullet,n_\bullet)= \sum_{i=1}^{t} n_i ( r \dim V_i - r_i \dim V)
= \sum_{i=1}^{t+1} \frac{\Gamma_i}{\dim V} ( r^i \dim V - r\dim
V^i)
\end{equation}
 (recall $n_i=\frac{\Gamma_{i+1}-\Gamma_i}{p}$), where
$r_{i}=\rk E_{i}$, $r^{i}=\rk E^{i}$, $E_{i}$ is the sheaf
generated by evaluation of the sections of $V_{i}$ and
$E^{i}=E_{i}/E_{i-1}$.

By the Hilbert-Mumford criterion (c.f. \cite[Theorem 2.1]{GIT} or \cite[Theorem 4.9]{Ne}), a point
$$
\overline{Q}\in \PP(\Hom(\wedge^r V , A) )
$$
is GIT (semi)stable if and only if for all weighted filtrations
$$
\mu(\overline{Q},V_\bullet,n_\bullet)(\leq) 0\; .
$$
Using the previous calculation, this can be stated as follows:

\begin{lem}
\label{HMcriterion}
A point $\overline{Q}$ is GIT (semi)stable if
for all weighted filtrations $(V_\bullet,n_\bullet)$
$$\sum_{i=1}^{t}  n_i ( r \dim V_i - r_i \dim V)(\leq) 0\; .$$
\end{lem}

The following theorem follows from \cite{Gi,Ma}

\begin{thm}
\label{equiv} Let $E$ be a sheaf. There exists an integer
$m_{0}(E)$ such that, for $m>m_{0}(E)$, the associated point
$\overline{Q}$ is GIT semistable if and only if the sheaf is
semistable.
\end{thm}

Let $E$ be an unstable sheaf. We choose an integer $m_{0}$ larger
than $m_0(E)$ and larger than
the integer used in Gieseker's construction of the moduli space.

Through Geometric Invariant Theory, stability of a point in the parameter space can be checked by $1$-parameter
subgroups (c.f. Hilbert-Mumford criterion, Proposition \ref{HMcriterion}): a point is unstable if there exists
any $1$-PS which makes some quantity positive. It is a natural question to ask if there exists a best way of
destabilizing a GIT unstable point, i.e. a best $1$-PS which gives maximum in the quantity we referred to.
Kempf explores this idea in \cite{Ke} and answers positively the
question, finding that there exists a special class
of $1$-parameter subgroup which moves most rapidly toward the origin.

We have seen that giving a weighted filtration, i.e. a filtration of vector subspaces $0\subset V_{1}\subset
\cdots \subset V_{t+1}= V$ and rational numbers $n_{1},\cdots,n_{t}>0$, is equivalent to giving a parabolic
subgroup with weights, which determines uniquely the vector $\Gamma$ of a $1$-PS and two of these $1$-PS are
conjugated by the parabolic and come from the same weighted filtration. We define the following function
$$
\nu(V_{\bullet},n_{\bullet})=\frac{\sum_{i=1}^{t}  n_{i} (r\dim V_{i}-r_{i}\dim V)}
{\sqrt{\sum_{i=1}^{t+1} {\dim V^{i}_{}} \Gamma_{i}^{2}}}\; ,
$$
which we call \emph{Kempf function}. Note that
$\nu(V_{\bullet},n_{\bullet})=\nu(V_{\bullet},\alpha
n_{\bullet})$, for every $\alpha>0$, hence by multiplying each
$n_{i}$ by the same scalar $\alpha$, which we call
\emph{rescaling the weights}, we get another 1-PS
but the same value for the Kempf function.

Note that this function corresponds to the one given in
\cite[Theorem 2.2]{Ke}. The numerator of both functions coincide
with the calculation of the minimal relevant weight by
Hilbert-Mumford criterion for GIT stability (c.f. \eqref{mrw}), and the denominator is the norm $||\Gamma ||$
of the vector
$$
\Gamma=(\overbrace{\Gamma_1,\ldots,\Gamma_1}^{\dim V^1},\overbrace{\Gamma_2,\ldots,\Gamma_2}^{\dim V^2},\ldots,\overbrace{\Gamma_{t+1},
\ldots,\Gamma_{t+1}}^{\dim V^{t+1}})
$$
as it is defined in \cite{Ke} as the Killing length of $\Gamma$.
Recall that for a simple group $G$ (as it is the case of
$G=\SL(V)$) every bilinear symmetric invariant form is a multiple
of the Killing form, and the norm $||\Gamma ||$ verifies these
properties.

We take the GIT quotient by the group $G=\SL(V)$, for which 
\cite[Theorem 2.2]{Ke}
 states that whenever there exists any $\Gamma$ giving
a positive value for the numerator of the function (i.e. whenever
there exists a 1-PS whose minimal relevant weight is positive,
which is equivalent to the sheaf $E$ to be unstable), there
exists a unique parabolic subgroup containing a unique
$1$-parameter subgroup in each maximal torus, giving maximum in
the Kempf function i.e., there exists a unique weighted filtration
for which the Kempf function achieves a maximum. Note that we
divide by the norm in the Kempf function to have
$\nu(V_{\bullet},n_{\bullet})=\nu(V_{\bullet},\alpha n_{\bullet}),
\forall \alpha>0$, hence a well defined maximal weighted
filtration for the function is defined up to rescaling, i.e., up
to multiplying every weight by the same positive scalar. 

Therefore, \cite[Theorem 2.2]{Ke} rewritten in our case asserts
the following:

\begin{thm}[Kempf]
\label{Kempf} There exists a unique weighted filtration (c.f. (\ref{filtV})) up to multiplication by a scalar, called \emph{Kempf filtration of
V}, such that the Kempf function $\nu(V_{\bullet},n_{\bullet})$ achieves the maximum among all filtrations and
weights $n_{i}>0$.
\end{thm}

We will construct a filtration by subsheaves of $E$
(which we will call Kempf filtration of $E$)
out of the Kempf
filtration of $V$.
Then we will relate the filtration given by Kempf
with the filtration constructed by Harder and Narasimhan to conclude
that both filtrations are the same.

\section{Convex cones}

In this subsection we define the machinery which will serve us in
the following. Endow $\mathbb{R}^{t+1}$ with an inner product
$(\cdot,\cdot)$ defined by a diagonal matrix
 $$
 \left(
 \begin{array}{ccc}
 b^1 & & 0 \\
  & \ddots & \\
 0 & & b^{t+1}\\
 \end{array}
 \right)
 $$
where $b^i$ are positive integers. Let
$$
\mathcal{C}= \big\{ x\in \mathbb{R}^{t+1} : x_1<x_2<\cdots <x_{t+1}
\big\}\; ,
$$
$$
\overline{\mathcal{C}}= \big\{ x\in \mathbb{R}^{t+1} :
x_1\leq x_2\leq \cdots \leq x_{t+1} \big\}\; ,
$$
and let $v= (v_1,\cdots,v_{t+1})\in \mathbb{R}^{t+1}-\{0\}$ verifying
$\sum_{i=0}^{t+1} v_{i}b^{i}=0$. Define the function
\begin{eqnarray*}
  \label{eq:mu}
\nu_{v}:\overline{\mathcal{C}}-\{0\} & \to & \mathbb{R}\\
\Gamma & \mapsto & \nu_{v}(\Gamma)=\frac{(\Gamma,v)}{||\Gamma||}
\end{eqnarray*}
and note that  $\nu_{v}(\Gamma)=||v||\cdot \cos\beta(\Gamma,v)$,
where $\beta(\Gamma, v)$ is the angle between $\Gamma$ and $v$.
Then, the function $\nu_{v}(\Gamma)$ does not depend on the norm
of $\Gamma$ and takes the same value on every point of the ray
spanned by each $\Gamma$.

Assume that there exists $\Gamma\in\overline{\mathcal{C}}$ with
$\nu_v(\Gamma)>0$. In that case, we want to find a vector
$\Gamma\in \overline{\mathcal{C}}$ which maximizes the function
defined before.

Let $w^i=-b^iv_i$, $w_0=0$, $w_i=w^1+\cdots+w^i$, $b_0=0$, and
$b_i=b^1+\cdots+b^i$. Note that $w_{t+1}=0$, by construction. We
draw a graph joining the points with coordinates $(b_i,w_i)$. Note
that this graph has $t+1$ segments, each segment has slope $-v_i$
and width $b^i$. This is the graph drawn with a thin line in the
figure. Now draw the convex envelope of this graph (thick line in
the figure), whose coordinates we denote by
$(b_i,\widetilde{w_i})$, and let us define
$\Gamma_{i}=-\frac{\widetilde{w_{i}}-\widetilde{w_{i-1}}}{b^i}$.
In other words, the quantities $-\Gamma_i$ are the slopes of the
convex envelope graph. We call the vector defined in this way
$\Gamma_{v}$. Note that the vector
$\Gamma_{v}=(\Gamma_1,\cdots,\Gamma_{t+1})$ belongs to
$\overline{\mathcal{C}}$ by construction and $\Gamma_{v}\neq 0$.

\setlength{\unitlength}{1cm}
\begin{picture}(13,7)(-1,-1)
\thicklines

\put(6,-0.3){\makebox(0,0)[c]{$b_i$}}
\put(-0.4,3){\makebox(0,0)[c]{\rotatebox{90}{$w_i,\widetilde{w_i}$}}}
\put(0,0){\makebox(0,0){$\circ$}}
\put(0,0){\line(1,1){4}}
\put(2,2){\makebox(0,0){$\circ$}}
\put(3,3){\makebox(0,0){$\circ$}}
\put(4,4){\makebox(0,0){$\circ$}}
\put(4,4){\line(4,1){4}}
\put(5,4.25){\makebox(0,0){$\circ$}}
\put(8,5){\makebox(0,0){$\circ$}}
\put(8,5){\line(1,-1){2}}
\put(9.5,3.5){\makebox(0,0){$\circ$}}
\put(10,3){\makebox(0,0){$\circ$}}
\put(10,3){\line(1,-3){1}}
\put(11,0){\makebox(0,0){$\circ$}}
\thinlines

\put(0,0){\line(1,0){11.5}}
\put(0,0){\line(0,1){5.5}}

\put(0,0){\makebox(0,0){$\circ$}}
\put(0,0){\line(2,1){2}}
\put(2,1){\makebox(0,0){$\circ$}}
\put(2,1){\line(1,1){1}}
\put(3,2){\makebox(0,0){$\circ$}}
\put(3,2){\line(1,2){1}}
\put(4,4){\line(1,-1){1}}
\put(5,3){\makebox(0,0){$\circ$}}
\put(5,3){\line(3,2){3}}
\put(9,4){\makebox(0,0){$\circ$}}
\put(9,4){\line(1,-4){0.5}}
\put(9.5,2){\makebox(0,0){$\circ$}}
\put(9.5,2){\line(1,2){0.5}}
\put(2,0.6){\makebox(0,0)[l]{$(b_1,w_1)$}}
\put(1.8,2.2){\makebox(0,0)[r]{$(b_1,\widetilde{w_1})$}}
\put(3,1.6){\makebox(0,0)[l]{$(b_2,w_2)$}}
\put(2.8,3.2){\makebox(0,0)[r]{$(b_2,\widetilde{w_2})$}}
\put(3.8,4.2){\makebox(0,0)[r]{$(b_3,\widetilde{w_3}=w_3)$}}
\put(5,4.7){\makebox(0,0)[c]{$(b_4,\widetilde{w_4})$}}
\put(5,2.6){\makebox(0,0)[c]{$(b_4,w_4)$}}
\put(8,5.4){\makebox(0,0)[c]{$(b_5,\widetilde{w_5}=w_5)$}}
\end{picture}

\begin{rem}
\label{rmk}
If $\widetilde{w_i}> w_i$, then $\Gamma_{i+1}=\Gamma_i$.
\end{rem}

\begin{thm}
\label{max}
The vector $\Gamma_{v}$ defined
in this way gives a maximum for the function $\nu_v$ on its domain.
\end{thm}

\begin{rem}
As one referee pointed out, this figure is very similar 
to \cite[Figure 2.1]{Br}, where for each
subsheaf $F$ of $E$ a point with coordinates $(\rk F,\deg F)$ is plotted,
and the Harder-Narasimhan polygon is the convex envelope.
Bruasse shows that the convex envelope achieves the maximum of certain
energy function coming from gauge theory. In our case, since we
are using Gieseker semistabiliby, instead of
the degree of $F$ we will plot $P_F(m)$, and the convex envelope
will maximize the function $\nu_v$ coming from GIT.
\end{rem}

Before proving the theorem we need some lemmas.

\begin{lem}
\label{close}
Let $v= (v_1,\cdots,v_{t+1})\in \mathbb{R}^{t+1}-\{0\}$ verifying
$\sum_{i=0}^{t+1} v_{i}b^{i}=0$.
Let $\Gamma$ be the point in $\overline{\mathcal{C}}$ which is closest
to $v$. Then $\Gamma$ achieves the maximum of $\nu_v$.
\end{lem}

\begin{pr}
For any $\alpha\in \mathbb{R}^{>0}$, the vector $\alpha \Gamma$ is also in
$\overline{\mathcal{C}}$, so in particular $\Gamma$ is the closest point in the
line $\alpha \Gamma$ to $v$. Therefore, the point
$\Gamma$ is the orthogonal
projection of $v$ into the line $\alpha \Gamma$, and the distance
is
\begin{equation}
  ||v|| \sin \beta(v,\Gamma)\label{eq:sin}
\end{equation}
where $\beta(\Gamma,v)$ is the angle between $\Gamma$ and $v$. But
a vector $\Gamma\in \overline{\mathcal{C}}$ minimizes
\eqref{eq:sin} if and only if it maximizes
$$
||v|| \cos \beta(\Gamma,v) = \frac{(\Gamma,v)}{||\Gamma||}\; ,
$$
so the lemma is proved.
\end{pr}

We say that an affine hyperplane in $\mathbb{R}^{t+1}$ separates a
point $v$ from $\mathcal{C}$ if $v$ is on one side of the
hyperplane and all the points of $\mathcal{C}$ are on
the other side of the hyperplane.

\begin{lem}
\label{hyperplane} Let $v\notin \overline{\mathcal{C}}$. A point
$\Gamma\in \overline{\mathcal{C}}-\{0\}$ gives minimum distance to $v$
if and only if the hyperplane $\Gamma+(v-\Gamma)^\perp$ separates
$v$ from $\mathcal{C}$.
\end{lem}

\begin{pr}
$\Rightarrow$) Let $\Gamma\in \overline{\mathcal{C}}$ and assume
there is a point $w\in \mathcal{C}$ on the same side of the
hyperplane as $v$. The segment going from $\Gamma$ to $w$ is in
$\overline{\mathcal{C}}$ (by convexity of
$\overline{\mathcal{C}}$), but there are points in this segment
(near $\Gamma$), which are closer to $v$ than $\Gamma$.

$\Leftarrow)$ Let $\Gamma$ be a point in $\overline{\mathcal{C}}$
such that $\Gamma+(v-\Gamma)^\perp$ separates $v$ from
$\mathcal{C}$. Let $w \in\overline{\mathcal{C}}$ be another point.
Let $w'$ be the intersection of the hyperplane and the segment
which goes from $w$ to $v$. Since the hyperplane separates
$\mathcal{C}$ from $v$, either $w'=w$ or $w'$ is in the interior
of the segment. Therefore
$$
d(w,v) \geq d(w',v) \geq d(\Gamma,v)
$$
where the last inequality follows from the fact that $\Gamma$
is the orthogonal projection of $v$ to the hyperplane.
\end{pr}

We thank F. Presas for suggesting this lemma, which is the key to prove Theorem \ref{max}.

\begin{pr}[Proof of the theorem]
Let $\Gamma_{v}=(\Gamma_{1},...,\Gamma_{t+1})$ be the vector in
the hypothesis of the theorem. If $v\in \overline{\mathcal{C}}$,
then $\Gamma_{v}=v$, and use Lemma \ref{close} to conclude. If
$v\notin \overline{\mathcal{C}}$, by Lemmas \ref{close} and
\ref{hyperplane}, it is enough to check that the hyperplane
$\Gamma_{v} + (v-\Gamma_{v})^\perp$ separates $v$ from
$\mathcal{C}$.

Let $\Gamma_{v}+\epsilon\in \mathcal{C}$. The condition that
$\Gamma_{v}+\epsilon$ belongs to $\mathcal{C}$ means that
\begin{equation}
  \label{eq:leq}
\epsilon_i-\epsilon_{i+1} < \Gamma_{i+1}-\Gamma_{i}
\end{equation}
The hyperplane separates $v$ from $\mathcal{C}$ if and only if
$(v-\Gamma_{v},\epsilon)<0$ for all such $\epsilon$. Therefore we
calculate (using the convention $\widetilde{w_0}=0$, $w_0=0$, and
$\widetilde{w_{t+1}}=w_{t+1}=0$)
$$
(v-\Gamma_{v},\epsilon) = \sum_{i=1}^{t+1}
b^i(v_i-\Gamma_i)\epsilon_i = \sum_{i=1}^{t+1}
(-w^i+(\widetilde{w_i}-\widetilde{w_{i-1}}))\epsilon_i =
$$
$$
=\sum_{i=1}^{t+1} \big( (\widetilde{w_i}-\widetilde{w_{i-1}})-(w_i-w_{i-1})\big)\epsilon_i
= \sum_{i=1}^{t+1} (\widetilde{w_{i}}- w_i)(\epsilon_i-\epsilon_{i+1})\; .
$$
If $\widetilde{w_i}=w_i$, then the corresponding summand is zero.
On the other hand, if $\widetilde{w_i}> w_i$, then
$\Gamma_{i+1}=\Gamma_i$ (c.f. Remark \ref{rmk}), and \eqref{eq:leq}
implies $\epsilon_i-\epsilon_{i+1}<0$. In any case, the summands
are always non-positive, and there is at least one which is
negative (because $v\notin \overline{\mathcal{C}}$ and then $v\neq
\Gamma_{v}$ and $\widetilde{w_i}> w_i$ for at least one $i$). Hence
$$(v-\Gamma_{v},\epsilon)<0\; .$$
\end{pr}

Therefore, the function $\nu_{v}(\Gamma)$ achieves its maximum for
the value $\Gamma_{v}\in \overline{\mathcal{C}}-\{0\}$ (or any other
point on the ray $\alpha\Gamma_{v}$) defined as the convex envelope
of the graph associated to $v$.

\section{Properties of the Kempf filtration}

Let $E$ be an unstable torsion-free sheaf over $X$ of Hilbert polynomial $P$. Let $m$ be an integer, $m\geq
m_{0}$, and let $V$ be a vector space of dimension $P(m)=h^{0}(E(m))$ (recall that $m_0$ was defined after
Theorem \ref{equiv}). We fix an isomorphism $V\simeq H^{0}(E(m))$ and let $0\subset V_{1}\subset \cdots \subset
V_{t+1}= V$ be the filtration of vector spaces given by Theorem \ref{Kempf}. Recall that it is called the
\emph{Kempf filtration of V}. For each index $i$, let $E^{m}_{i}\subset E$ be the subsheaf generated by $V_{i}$
under the evaluation map. We call this filtration
$$
0\subseteq E^{m}_1 \subseteq E^{m}_2 \subseteq \;\cdots\; \subseteq E^{m}_t \subseteq E^{m}_{t+1}=E\; ,
$$
the \emph{$m$-Kempf filtration of E}.

Recall that given a filtration of $V$, and the corresponding filtration of $E$ obtained by evaluating, we call $V^{i}=V_{i}/V_{i-1}$, 
$E^{i}=E_{i}/E_{i-1}$, and $r_{i}=\rk E_{i}$, $r^{i}=\rk E^{i}$. 

\begin{dfn}
\label{graph} Let $m\geq m_{0}$ and let $0\subset V_{1}\subset \cdots \subset V_{t+1}=V$ be a filtration of vector
spaces of $V$. Let
$$v_{m,i}=m^{n+1}\cdot \frac{1}{\dim V^{i}\dim V}\big[r^{i}\dim V-r\dim V^{i}\big]\; ,$$
$$b_{m}^{i}=\dfrac{1}{m^{n}}\dim V^{i}>0\; ,$$
$$w_{m}^{i}=-b_{m}^{i}\cdot v_{m,i}=m\cdot \frac{1}{\dim V}\big[r\dim V^{i}-r^{i}\dim V\big]\; .$$
Also let
$$b_{m,i}=b_{m}^{1}+\ldots +b_{m}^{i}=\dfrac{1}{m^{n}}\dim V_{i}\; ,$$
$$w_{m,i}=w_{m}^{1}+\ldots +w_{m}^{i}=m\cdot \frac{1}{\dim V}\big[r\dim V_{i}-r_{i}\dim V\big]\; .$$
We call the graph defined by points $(b_{m,i},w_{m,i})$ the \emph{graph associated to the filtration} $V_{\bullet}\subset V$.
\end{dfn}

Now we can identify the Kempf function, i.e. the function in Theorem \ref{Kempf}
$$\nu(V_{\bullet},n_{\bullet})=\frac{\sum_{i=1}^{t} n_{i}(r\dim V_{i}-r_{i}\dim V)}{\sqrt{\sum_{i=1}^{t+1}
{\dim V^{i}} \Gamma_{i}^{2}}}\; ,$$ with the function in Theorem
\ref{max} up to a factor which is a power of $m$, by defining
$v_{m,i}$, the coordinates of vector $v_{m}$, and $b_{m}^{i}$, the
eigenvalues of the scalar product, as in Definition \ref{graph}.
Note that $-v_{m,i}$ are the slopes of the graph associated to the
filtration $V_{\bullet}\subset V$. Here the coordinates
$\Gamma_{i}$ are the same as in the $1$-PS defined by $n_{i}$.
Also note that $\sum_{i=1}^{t+1}v_{m,i}b^{i}_{m}=0$. Then, an easy
calculation shows that

\begin{prop}
\label{identification}
For every integer $m$, the following equality holds
$$
\nu(V_{\bullet},n_{\bullet})=m^{(-\frac{n}{2}-1)}\cdot
\nu_{v_{m}}(\Gamma)
$$
between the Kempf function in Theorem \ref{Kempf} and the function in Theorem \ref{max}.
\end{prop}

In the following, we will omit the subindex $m$ for the numbers
$v_{m,i}$, $b_{m,i}$, $w_{m,i}$ in the definition of the graph
associated to a filtration of vector spaces, where it is clear
from the context. Hence, given $V\simeq H^{0}(E(m))$ we will refer
to a filtration $V_{\bullet}\subset V$ and a vector
$v=(v_{1},\ldots,v_{t+1})$ as the vector of the graph associated
to the filtration.

\begin{rem}
We introduce the factor $m^{n+1}$ in Definition \ref{graph} for
convenience, so that $v_{m,i}$ and $b_{m}^{i}$ have order zero on
$m$, because $\dim V=P(m)$ appears in their expressions. Then, the
size of the graph does not change when $m$ grows.
\end{rem}

\begin{lem}
\label{lemmaA}
Let $0\subset V_{1}\subset \cdots \subset V_{t+1}= V$ be the Kempf filtration of $V$ (cf. Theorem \ref{Kempf}).
Let $v=(v_{1},...,v_{t+1})$ be the vector of the graph associated to this filtration by Definition \ref{graph}.
Then
$$
v_{1}<v_{2}<\ldots<v_{t}<v_{t+1}\; ,
$$
i.e., the graph is convex.
\end{lem}

\begin{pr}
By Theorem \ref{Kempf} the maximum of $\nu$ among all
filtrations $V_{\bullet}\subset V$ and weights $n_{i}>0,\forall i$
is achieved by a unique weighted filtration
$(V_{\bullet},n_{\bullet})$, $n_{i}>0,\forall i$, up to rescaling. Let
$V_{\bullet}\subset V$ be this filtration, and allow $n_{i}$ to
vary. By Proposition \ref{identification}
$\nu$ is equal to $\nu_v$ up to a constant factor.
By Theorem \ref{max}, $\nu_v$
achieves the maximum on $\Gamma_{v}$. The vector
$\Gamma_{v}$ corresponds to the weights $n_{i}$ given by Theorem
\ref{Kempf}. Summing up, if $V_\bullet\subset V$ is Kempf filtration
of $V$, then the vector
$\Gamma_{v}=(\Gamma_{1},\ldots,\Gamma_{t+1})$ has
$\Gamma_{i}<\Gamma_{i+1},\forall i$.

Assume that, for the Kempf filtration of $V$, there exists some $i$ such
that $v_{i}\geq v_{i+1}$. Then $v\notin \mathcal{C}$ and, by Lemma
\ref{close}, $\Gamma_{v}\in \overline{\mathcal{C}}\backslash
\mathcal{C}$, which means that there exists some $j$ with
$\Gamma_{j}=\Gamma_{j+1}$, but we have just seen that this is
impossible.
\end{pr}

\begin{lem}
\label{lemmaB} Let $0\subset V_{1}\subset \cdots \subset V_{t+1}= V$ be the Kempf filtration of $V$ (cf. Theorem
\ref{Kempf}). Let $W$ be a vector space with $V_{i}\subset W\subset V_{i+1}$ and consider the new filtration
$V'_{\bullet}\subset V$
\begin{equation}
    \begin{array}{ccccccccccccccccc}
    0 & \subset & V'_{1} & \subset & \cdots & \subset & V'_{i} & \subset & V'_{i+1} & \subset & V'_{i+2} & \subset & \cdots & \subset & V'_{t+2} & =& V\\
    || & & || & & & & || & & || & & || & & & & & & ||  \\
    0 & \subset & V_{1}& \subset & \cdots & \subset & V_{i} & \subset & W & \subset & V_{i+1} & \subset & \cdots & \subset & V_{t+1} & = & V
    \end{array}
\end{equation}
Then, $v'_{i+1}\geq v_{i+1}$. We say that the Kempf filtration is
the convex envelope of every refinement.
\end{lem}

\begin{pr}
The graph associated to $V'_{\bullet}\subset V$ has one more point
than the graph associated to $V_{\bullet}\subset V$, hence it is a
refinement of the graph associated to Kempf filtration of $V$.
Therefore the convex envelope of the graph associated to $v'$ has
to be equal to the graph associated to $v$, and this happens only
when the extra point associated to $W$ is not above the graph
associated to $v$, which means that the slope $-v'_{i+1}$ has to
be less or equal than $-v_{i+1}$.
\end{pr}

Later on, we will check that, for $m$ large enough, the $m$-Kempf
filtration stabilizes in the sense $E^{m}_{i}=E^{m+l}_{i},\forall
i,\forall l>0$. The $m$-Kempf filtration for $m\gg 0$ will be
called the Kempf filtration of $E$, and the goal of this article
is to show that it coincides with the Harder-Narasimhan filtration
of $E$.

\begin{lem}[Simpson]
\cite[Corollary 1.7]{Si} or \cite[Lemma 2.2]{HL1}
\label{LePSi} Let $r>0$ be an integer. Then there exists a
constant $B$ with the following property: for every torsion free
sheaf $E$ with $0<\rk (E)\leq r$, we have
$$h^{0}(E)\leq \frac{1}{g^{n-1}n!}\big ((\rk(E)-1)([\mu_{max}(E)+B]_{+})^{n}+([\mu_{mim}(E)+B]_{+})^{n}\big)\; ,$$
where $g=\deg \mathcal{O}_{X}(1)$, $[x]_{+}=\max\{0,x\}$, and
$\mu_{max}(E)$ (respectively $\mu_{min}(E)$) is the maximum (resp.
minimum) slope of the Mumford-semistable factors of the
Harder-Narasimhan filtration of $E$.
\end{lem}

\begin{rem}
Recall that the Harder-Narasimhan filtration with Gieseker
stability is a refinement of the one with Mumford stability, with
the inequalities holding between polynomials of their leading
coefficients.
\end{rem}

We denote
$$P_{\mathcal{O}_{X}}(m)=\frac{\alpha_{n}}{n!}m^{n}+\frac{\alpha_{n-1}}{(n-1)!}m^{n-1}+...
+\frac{\alpha_{1}}{1!}m+\frac{\alpha_{0}}{0!}$$ the Hilbert
polynomial of $\mathcal{O}_{X}$, then $\alpha_{n}=g$. Let
$$P(m)=\frac{rg}{n!}m^{n}+\frac{d+r\alpha_{n-1}}{(n-1)!}m^{n-1}+...$$
be the Hilbert polynomial of the sheaf $E$, where $d$ is the degree and $r$ is the rank. Let us call
$A=d+r\alpha_{n-1}$, so
$$P(m)=\frac{rg}{n!}m^{n}+\frac{A}{(n-1)!}m^{n-1}+...$$
Let us define
\begin{equation}
\label{C_constant}
C=\max\{r|\mu _{\max }(E)|+\frac{d}{r}+r|B|+|A|+1\;,\;1\},
\end{equation}
a positive constant.

\begin{prop}
\label{boundedness} Given an integer $m$ and a vector space $V\simeq H^{0}(E(m))$, we have the Kempf filtration
$V_{\bullet}\subset V\simeq H^{0}(E(m))$ and, by evaluation, the $m$-Kempf filtration $E_{\bullet}^{m}\subseteq
E$. There exists an integer $m_{2}$ such that for $m\geq m_{2}$, each term in the $m$-Kempf filtration of $E$
has slope $\mu(E^{m}_{i})\geq \dfrac{d}{r}-C$.
\end{prop}

\begin{pr}
Choose an $m_{1}\geq m_{0}$ such that for $m\geq m_{1}$
$$[\mu_{max}(E)+gm+B]_{+}=\mu_{max}(E)+gm+B$$
and
$$[\frac{d}{r}-C+gm+B]_{+}=\frac{d}{r}-C+gm+B\; .$$
Now let $m\geq m_{1}$ and let
$$0 \subseteq E^{m}_1 \subseteq E^{m}_2 \subseteq \;\cdots\; \subseteq E^{m}_t \subseteq E^{m}_{t+1}=E$$
be the $m$-Kempf filtration.

Suppose we have a term in the filtration $E^{m}_{i}\subseteq E$, of rank $r_{i}$ and degree $d_{i}$, such that $\mu
(E^{m}_{i})<\frac{d}{r}-C$. The subsheaf $E_{i}^{m}(m)\subseteq E(m)$ satisfies the estimate in Lemma
\ref{LePSi},
$$h^{0}(E^{m}_{i}(m))\leq \frac{1}{g^{n-1}n!}\big ((r_{i}-1)([\mu_{max}(E^{m}_{i})+gm+B]_{+})^{n}+([\mu_{min}(E^{m}_{i})+gm+B]_{+})^{n}\big)\; ,$$
where $\mu_{max}(E^{m}_{i}(m))=\mu_{max}(E^{m}_{i})+gm$ and similarly for $\mu_{min}$.

Note that $\mu_{max}(E^{m}_{i})\leq \mu_{max}(E)$ and
$\mu_{min}(E^{m}_{i})\leq \mu(E^{m}_{i})<\frac{d}{r}-C$, so
$$h^{0}(E^{m}_{i}(m))\leq \frac{1}{g^{n-1}n!}\big ((r_{i}-1)([\mu_{max}(E)+gm+B]_{+})^{n}+([\frac{d}{r}-C+gm+B]_{+})^{n}\big)\; ,$$
and, by choice of $m$,
$$h^{0}(E^{m}_{i}(m))\leq \frac{1}{g^{n-1}n!}\big ((r_{i}-1)(\mu_{max}(E)+gm+B)^{n}+(\frac{d}{r}-C+gm+B)^{n}\big)=G(m)\; ,$$
where
$$G(m)=\frac{1}{g^{n-1}n!}\big
[r_{i}g^{n}m^{n}+ng^{n-1}\big((r_{i}-1)\mu_{max}(E)+\frac{d}{r}-C+r_{i}B\big)m^{n-1}+\cdots
\big ]\; .$$

Recall that, by Definition \ref{graph}, to such filtration we
associate a graph with heights, for each $j$,
$$w_{j}=w^{1}+\ldots +w^{j}=m\cdot \frac{1}{\dim V}\big[r\dim V_{j}-r_{j}\dim V\big]\; .$$
To reach a contradiction, it is enough to show that $w_{i}<0$. In
that case, the graph has to be convex by Lemma \ref{lemmaA}. If
$w_{i}<0$ there is a $j<i$ such that $-v_{j}<0$, because the graph
starts on the origin. Hence, the rest of the slopes of the graph
are negative, $-v_{k}<0$, $k\geq i$, because the slopes have to be
decreasing. Then $w_{i}>w_{i+1}>\ldots w_{t+1}$, and $w_{t+1}<0$.
But it is
$$w_{t+1}=m\cdot \frac{1}{\dim V}\big[r\dim V_{t+1}-r_{t+1}\dim V\big]=0\; ,$$
because $r_{t+1}=r$ and $V_{t+1}=V$, then the contradiction.

Let us show that $w_{i}<0$. Since $E^{m}_{i}(m)$ is generated by
$V_{i}$ under the evaluation map, it is $\dim V_{i}\leq
h^{0}(E^{m}_{i}(m))$, hence
$$
w_{i}=\frac{m}{\dim V}\big[r\dim V_{i}-r_{i}\dim V\big]\leq
$$
$$
\leq \frac{m}{P(m)}\big[rh^{0}(E^{m}_{i}(m))-r_{i}P(m)\big]\leq
\frac{m}{P(m)}\big[rG(m)-r_{i}P(m)\big]\; .
$$

Hence, $w_{i}<0$ is equivalent to
$$\Psi(m)=rG(m)-r_{i}P_{E}(m)<0\; ,$$
where $\Psi(m)=\xi_{n}m^{n}+\xi_{n-1}m^{n-1}+\cdots +\xi_{1}m+\xi_{0}$ is an $n$-order polynomial. Let us
calculate the $n^{th}$-coefficient:
$$\xi_{n}=(rG(m)-r_{i}P(m))_{n}=r\frac{r_{i}g}{n!}-r_{i}\frac{rg}{n!}=0\; .$$
Then, $\Psi(m)$ has no coefficient in order $n$. Let us calculate the $(n-1)^{st}$-coefficient:
$$\xi_{n-1}=(rG(m)-r_{i}P(m))_{n-1}=(rG_{n-1}-r_{i}\frac{A}{(n-1)!})$$
where $G_{n-1}$ is the $(n-1)^{st}$-coefficient of the polynomial $G(m)$,
$$G_{n-1}=\frac{1}{g^{n-1}n!}ng^{n-1}((r_{i}-1)\mu_{max}(E)+\frac{d}{r}-C+r_{i}B)=$$
$$\frac{1}{(n-1)!}((r_{i}-1)\mu_{max}(E)+\frac{d}{r}-C+r_{i}B)\leq$$
$$\frac{1}{(n-1)!}((r_{i}-1)|\mu_{max}(E)|+\frac{d}{r}-C+r_{i}|B|)\leq$$
$$\frac{1}{(n-1)!}(r|\mu_{max}(E)|+\frac{d}{r}-C+r|B|)<\frac{-|A|}{(n-1)!}\; ,$$
last inequality coming from the definition of $C$ in \eqref{C_constant}.
Then
$$\xi_{n-1}<r\big(\frac{-|A|}{(n-1)!}\big)-r_{i}\frac{A}{(n-1)!}=
\frac{-r|A|-r_{i}A}{(n-1)!}<0$$ because $-r|A|-r_{i}A<0$.

Therefore $\Psi(m)=\xi_{n-1}m^{n-1}+\cdots +\xi_{1}m+\xi_{0}$
with $\xi_{n-1}<0$, so there exists $m_{2}\geq m_{1}$ such that for $m\geq
m_{2}$ we will have $\Psi(m)<0$ and
$w_{i}<0$, then the contradiction.
\end{pr}

\begin{prop}
\label{regular}
There exists an integer $m_{3}$ such that for $m\geq m_{3}$ the sheaves $E^{m}_{i}$ and $E^{m,i}=E^m_i/E^m_{i-1}$ are $m_{3}$-regular. In particular their higher cohomology groups vanish and they are generated by global sections.
\end{prop}
\begin{pr}
Note that $\mu(E^{m}_{i})\leq \mu_{\max}(E)$. Then, although $E^{m}_{i}$ depends on $m$, its slope is bounded above and below by numbers which do not depend on $m$, (cf. Proposition \ref{boundedness}) and furthermore it is a subsheaf of $E$. Hence, the set of possible isomorphism classes for $E^{m}_{i}$ is bounded.
Apply Serre vanishing theorem choosing $m_{3}\geq m_{2}$.
\end{pr}


\begin{prop}
Let $m\geq m_{3}$. For each term $E_{i}^{m}$ in the $m$-Kempf filtration, we have $\dim V_{i}=h^{0}(E_{i}^{m}(m))$, therefore $V_{i}\cong H^{0}(E_{i}^{m}(m))$.
\label{task}
\end{prop}

\begin{pr}
Let $V_{\bullet}\subset V$ be the Kempf filtration of $V$ (cf. Theorem \ref{Kempf}) and let
$E_{\bullet}^{m}\subseteq E$ be the $m$-Kempf filtration of $E$. We know that each $V_{i}$ generates the
subsheaf $E_{i}^{m}$, by definition, then we have the following diagram:

$$\begin{array}{ccccccccccc}
    0 & \subset & V_{1} & \subset & V_{2} & \subset & \cdots & \subset & V_{t+1} & = & V \\
     & & \cap & & \cap & & & & & & ||  \\
      &     & H^{0}(E_{1}^{m}(m)) & \subset & H^{0}(E_{2}^{m}(m)) & \subset & \cdots & \subset & H^{0}(E_{t+1}^{m}(m)) & = & H^{0}(E(m))
      \end{array}$$

Suppose there exists an index $i$ such that $V_{i}\neq H^{0}(E_{i}^{m}(m))$. Let $i$ be the index such that $V_{i}\neq H^{0}(E_{i}^{m}(m))$ and $\forall j>i$ it is $V_{j}=H^{0}(E_{j}^{m}(m))$. Then we have the diagram:
\begin{equation}
 \label{filtrationV}
    \begin{array}{ccccccccccccc}
    V_{i} & \subset & V_{i+1}\\
    \cap & & ||\\
    H^{0}(E_{i}^{m}(m)) & \subseteq & H^{0}(E_{i+1}^{m}(m))
    \end{array}
\end{equation}

Therefore $V_{i}\subsetneq H^{0}(E_{i}^{m}(m))\subsetneq V_{i+1}$ and we can consider a new filtration by adding the term $H^{0}(E_{i}^{m}(m))$:
\begin{equation}
\label{filtrationV'}
    \begin{array}{ccccccccccccc}
    V_{i} & \subset & H^{0}(E_{i}^{m}(m)) & \subset & V_{i+1}\\
    || & & || & & ||\\
    V'_{i} & &  V'_{i+1} & & V'_{i+2}
    \end{array}
\end{equation}

Note that we are in situation of Lemma \ref{lemmaB}, where
$W=H^{0}(E_{i}^{m}(m))$, filtration $V_{\bullet}$ is
\eqref{filtrationV} and filtration $V'_{\bullet}$ is
\eqref{filtrationV'}.

The graph associated to filtration $V_{\bullet}$, by Definition \ref{graph}, is given by the points
$$(b_{i},w_{i})=\big(\dfrac{\dim V_{i}}{m^{n}},\frac{m}{\dim V}(r\dim V_{i}-r_{i}\dim V)\big)\; ,$$
where the slopes of the graph are given by

$$-v_{i}=\frac{w^{i}}{b^{i}}=\frac{w_{i}-w_{i-1}}{b_{i}-b_{i-1}}=$$
$$\frac{m^{n+1}}{\dim V}\big(r-r^{i}\frac{\dim V}{\dim V^{i}}\big)\leq \frac{m^{n+1}}{\dim V}\cdot r:=R$$
and equality holds if and only if $r^{i}=0$.

Now, the new point which appears in the graph of the filtration $V'_{\bullet}$ is
$$Q=\big(\dfrac{h^{0}(E_{i}^{m}(m))}{m^{n}},\frac{m}{\dim V}(rh^{0}(E_{i}^{m}(m))-r_{i}\dim V)\big)\; .$$

Point $Q$ joins two new segments appearing in this new graph. The slope of the segment between $(b_{i},w_{i})$ and $Q$ is, by a
similar calculation,
$$-v'_{i+1}=\dfrac{m^{n+1}}{\dim V}\cdot r=R\; .$$

By Lemma \ref{lemmaA}, the graph is convex, hence
$v_{1}<v_{2}<\ldots<v_{t+1}$. As $E^{m}_{1}$ is a non zero
torsion-free sheaf, it has positive rank $r_{1}=r^{1}$ and so it
follows $v_{1}>-R$. On the other hand, by Lemma \ref{lemmaB},
$v'_{i+1}\geq v_{i+1}$. Hence
$$-R<v_{1}<v_{2}<\ldots<v_{i+1}\leq v'_{i+1}=-R\; ,$$
which is a contradiction.

Therefore, $\dim V_{i}=h^{0}(E_{i}^{m}(m))$, for every term in
the $m$-Kempf filtration.
\end{pr}

\begin{cor}
\label{rank}
For every term $E_{i}^{m}$ in the $m$-Kempf filtration, we have $r^{i}>0$, where $r^{i}=\rk E^{i,m}=\rk E_{i}^{m}/E_{i-1}^{m}$.
\end{cor}
\begin{pr}
We have seen that $r^{i}=0$ is equivalent to $-v_{i}=R$. Then, the result follows from Proposition \ref{task} because it is 
$r^{1}=r_{1}>0$ and $-R<v_{1}<v_{2}<\cdots<v_{t+1}$.
\end{pr}

\section{The $m$-Kempf filtration stabilizes}

In Proposition \ref{regular} we have seen that, for any $m\geq m_{3}$,
all the terms $E^{m}_{i}$ in the $m$-Kempf filtration of $E$ are
$m_{3}$-regular. Hence, $E^{m}_{i}(m_{3})$ is generated by the
subspace $H^{0}(E^{m}_{i}(m_{3}))$ of $H^{0}(E(m_{3}))$, and the
filtration of sheaves
$$
0\subset E^{m}_{1} \subset E^{m}_{2} \subset \cdots \subset E^{m}_{t_{m}} \subset E^{m}_{t_{m}+1}=E
$$
is the filtration associated to the filtration of vector spaces
$$
0\subset H^{0}(E^{m}_{1}(m_{3})) \subset H^{0}(E^{m}_{2}(m_{3})) \subset \cdots \subset
H^{0}(E^{m}_{t_{m}}(m_{3})) \subset H^{0}(E^{m}_{t_{m}+1}(m_{3}))=H^{0}(E(m_{3}))
$$

y the evaluation map. Note that the dimension of the vector space
$H^{0}(E(m_{3}))$ does not depend on $m$ and, by Corollary
\ref{rank}, the length $t_{m}+1$ of the $m$-Kempf filtration of
$E$ is at most equal to $r$, the rank of $E$, a bound which does
not also depend on $m$.

We call \emph{$m$-type} to the tuple of different Hilbert
polynomials appearing in the $m$-Kempf filtration of $E$
$$
(P_{1}^{m},\ldots,P_{t_{m}+1}^{m})\; ,
$$
where $P_{i}^{m}:=P_{E_{i}^{m}}$. Note that
$P^{i,m}:=P_{E^{m}_{i}/E^{m}_{i-1}}=P_{E^{m}_{i}}-P_{E^{m}_{i-1}}$,
so they are defined in terms of elements of each $m$-type.

\begin{prop}
\label{Pisfinite}
For all integers $m\geq m_{3}$, the set of
possible $m$-types
$$
\mathcal{P}=\big\{(P_{1}^{m},\ldots, P_{t_{m}+1}^{m})\big\}
$$
is finite.
\end{prop}
\begin{pr}
Once we fix $V\cong H^{0}(E(m_{3}))$ of dimension
$h^{0}(E(m_{3}))$ (which does not depend on $m$), all the possible
filtrations by vector subspaces are parametrized by a finite-type
scheme. Therefore the set of all possible $m$-Kempf filtrations of
$E$, for $m\geq m_{3}$, is bounded and $\mathcal{P}$ is finite.
\end{pr}

Recall that the vector $v$ can be recovered from the filtration
$V_{\bullet}\subset V$ and the vector $\Gamma$ from the weights
$n_{i}$. Then, given $m$, the $m$-Kempf filtration achieves the
maximum for the function $\nu(V_{\bullet},n_{\bullet})$, which is the
same, by Proposition \ref{identification}, as achieving the maximum
for the function
$$
\nu_v(\Gamma)=\frac{(\Gamma,v)}{||\Gamma||}\; ,
$$
among all filtrations $V_{\bullet}\subset V$ and vectors $\Gamma\in
\mathcal{C}-\{0\}$, where
$$
\mathcal{C}= \big\{ x\in \mathbb{R}^{t+1} : x_1<x_2<\cdots
<x_{t+1} \big\}\; .
$$

By Definition \ref{graph} we associate a graph to the $m$-Kempf
filtration, given by $v_{m}$. Recall that, by Lemma \ref{lemmaA},
the graph is convex, meaning $v_{m}\in \mathcal{C}$, which implies
$\Gamma_{v_{m}}=v_{m}$ by Lemma \ref{close}. Then, given $v_{m}$
associated to the $m$-Kempf filtration,
\begin{equation}
\label{maxvalue}
\max_{\Gamma\in \overline{\mathcal{C}}}\nu_{v_{m}}(\Gamma)=\nu_{v_{m}}(\Gamma_{v_{m}})=\frac{(\Gamma_{v_{m}},v_{m})}{||\Gamma_{v_{m}}||}=
\frac{(v_{m},v_{m})}{||v_{m}||}=||v_{m}||\; ,
\end{equation}
where recall that we defined (c.f. Definition \ref{graph})
$$
v_{m,i}=m^{n+1}\cdot \frac{1}{\dim V^{i}\dim V}\big[r^{i}\dim V-r\dim
V^{i}\big]
$$
$$
b_{m}^{i}=\frac{1}{m^{n}}\cdot \dim V^{i}
$$
and, thanks to Propositions \ref{regular} and \ref{task}, we can rewrite
$$
v_{m,i}=m^{n+1}\cdot \frac{1}{P^{i,m}(m)P(m)}\big[r^{i}P(m)-rP^{i,m}(m)\big]\; ,
$$
$$
b_{m}^{i}=\frac{1}{m^{n}}\cdot P^{i,m}(m)\; .
$$

Let
$$
v_{m,i}(l)=m^{n+1}\cdot \frac{1}{P^{i,m}(l)P(l)}\big[r^{i}P(l)-rP^{i,m}(l)\big]\; ,
$$
and let us define
$$
\Theta_{m}(l)=(\nu_{v_{m}(l)}(\Gamma_{v_{m}(l)}))^{2}=||v_{m}(l)||^{2}\;
,
$$
where the second equality follows by an argument similar to
\eqref{maxvalue}. Note that $\Theta_m(l)$ is a rational
function on $l$. Let
$$
\mathcal{A}=\{\Theta_{m}:m\geq m_{3}\}\; ,
$$
which is a finite set by Proposition \ref{Pisfinite}. We say that
$f_{1}\prec f_{2}$ for two rational functions, if the inequality

$f_{1}(l)<f_{2}(l)$ holds for $l\gg 0$, and let $K$ be the maximal
function in the finite set $\mathcal{A}$, with respect to the
defined ordering.

Note that the value $\Theta_{m}(m)$ is the square of the maximum
of Kempf's function $\nu_{v_{m}}(\Gamma)$, by \eqref{maxvalue},
achieved for the maximal filtration $V_{\bullet}\subset V\simeq
H^{0}(E(m))$ of vector spaces which gives the vector $v_{m}$. This
weighted filtration is the only one which gives the value
$\sqrt{\Theta_{m}(m)}$ for the Kempf function.

\begin{lem}
\label{uniquefunction} There exists an integer $m_{4}\geq m_{3}$
such that $\forall m\geq m_{4}$, $\Theta_{m}=K$.
\end{lem}
\begin{pr}
Choose $m_{4}$ such that $K(l) \geq \Theta_{m}(l)$, $\forall l\geq
m_{4}$ and every $\Theta_{m}\in \mathcal{A}$ with equality only
when $\Theta_{m}= K$, and let $m\geq m_{4}$. Given that the Kempf
function achieves the maximum over all possible filtrations and
weights (c.f. Theorem \ref{Kempf}), we have $\Theta_{m}(m)\geq
K(m)$, because $K$ is another rational function built with other
$m'$-type, i.e., other values for the polynomials appearing on the
rational function. Combining both inequalities we obtain
$\Theta_m(m)=K(m)$ for all $m\geq m_{4}$.
\end{pr}

\begin{prop}
\label{eventually} Let $l_{1}$ and $l_{2}$ be integers with
$l_{1}\geq l_{2} \geq m_{4}$. Then the $l_{1}$-Kempf filtration of
$E$ is equal to the $l_{2}$-Kempf filtration of $E$.
\end{prop}

\begin{pr}
By construction, the filtration
\begin{equation}
  \label{filt1}
  0\subset H^{0}(E^{l_{1}}_{1}(l_{1})) \subset H^{0}(E^{l_{1}}_{2}(l_{1})) \subset \cdots \subset H^{0}(E^{l_{1}}_{t_{1}}(l_{1}))
  \subset H^{0}(E^{l_{1}}_{t_{1}+1}(l_{1}))=H^{0}(E(l_{1}))
\end{equation}
is the $l_{1}$-Kempf filtration of $V\simeq H^{0}(E(l_{1}))$. Now
consider the filtration $V'_{\bullet}\subset V\simeq
H^{0}(E(l_{1}))$ defined as follows
\begin{equation}
\label{filt2} 0\subset H^{0}(E^{l_{2}}_{1}(l_{1})) \subset H^{0}(E^{l_{2}}_{2}(l_{1})) \subset \cdots \subset
H^{0}(E^{l_{2}}_{t_{2}}(l_{1})) \subset H^{0}(E^{l_{2}}_{t_{2}+1}(l_{1}))=H^{0}(E(l_{1})) \; .
\end{equation}
We have to prove that \eqref{filt2} is in fact the $l_{1}$-Kempf
filtration of $V\simeq H^{0}(E(l_{1}))$.

Since $l_{1},l_{2}\geq m_{4}$, by Lemma \ref{uniquefunction} we
have $\Theta_{l_{1}}=\Theta_{l_{2}}=K$. Then,
$\Theta_{l_{1}}(l_{1})=\Theta_{l_{2}}(l_{1})$ and, by uniqueness
of the Kempf filtration (c.f. Theorem \ref{Kempf}), the
filtrations
 \eqref{filt1} and \eqref{filt2} coincide. Since, in particular $l_{1},l_{2}\geq m_{3}$, $E_{i}^{l_{1}}$ and $E_{i}^{l_{2}}$
 are $l_{1}$-regular by Proposition \ref{regular}. Hence, $E_{i}^{l_{1}}(l_{1})$ and $E_{i}^{l_{2}}(l_{1})$ are generated by their global sections (c.f. Lemma \ref{mregularity})
 $H^{0}(E_{i}^{l_{1}}(l_{1}))=H^{0}(E_{i}^{l_{2}}(l_{1}))$, which are equal by the
 previous argument, therefore $E_{i}^{l_{1}}(l_{1})=E_{i}^{l_{2}}(l_{1})$. By tensoring with $\mathcal{O}_{X}(-l_{1})$, this implies that
 the filtrations $E^{l_{1}}_{\bullet} \subset E$
  and $E^{l_{2}}_{\bullet} \subset E$ coincide.
\end{pr}

\begin{dfn}
If $m\geq m_{4}$, the $m$-Kempf filtration of $E$ is called
\emph{the Kempf filtration of $E$},
$$0\subset E_{1} \subset E_{2} \subset \cdots \subset E_{t} \subset E_{t+1}=E\; .$$
\end{dfn}

\section{Kempf filtration is Harder-Narasimhan filtration}

Recall that the Kempf theorem (c.f. Theorem \ref{Kempf}) asserts
that given an integer $m$ and $V\simeq H^{0}(E(m))$, there exists a
unique weighted filtration of vector spaces $V_{\bullet}\subseteq V$
which gives maximum for the Kempf function
$$
\nu(V_{\bullet},n_{\bullet})=
\frac {\sum_{i=1}^{t+1} \frac{\Gamma_i}{\dim V} ( r^i \dim V - r\dim V^i)}
{\sqrt{\sum_{i=1}^{t+1} {\dim V^{i}} \Gamma_{i}^{2}}}\; .
$$
This filtration induces a filtration of sheaves, called the Kempf
filtration of $E$,
$$
0\subset E_{1} \subset E_{2} \subset \cdots \subset E_{t} \subset
E_{t+1}=E
$$
which is independent of $m$, for $m\geq m_{4}$, by Proposition
\ref{eventually}, hence it only depends on $E$. From now on, we assume
$m\geq m_4$.

In the previous sections, based on the fact we can rewrite the Kempf function as a certain scalar product
divided by a norm (c.f. Proposition \ref{identification}), we saw that Kempf filtration is encoded by a convex
graph (c.f. Lemma \ref{lemmaA}). We can express the data related to the filtration of vector spaces with the
data of the corresponding filtration of sheaves. Since $m\geq m_3$, the sheaves $E_{i}$ and $E^i$ are
$m$-regular $\forall i$, and
\begin{equation}
\label{substitutions}
    \begin{array}{c}
    \dim V_{i}=h^{0}(E_{i}(m))=P_{E_{i}}(m) =: P_i(m)\\
\dim V^{i}=h^{0}(E^{i}(m))=P_{E^{i}}(m) =: P^i(m)
    \end{array}
\end{equation}
(c.f. Proposition \ref{regular} and Proposition \ref{task}). Recall that the Kempf function is a rational
function on $m$, with order $m^{-\frac{n}{2}-1}$ at zero (c.f. Proposition \ref{identification}). Then we
consider the function $\nu$, where
$$
\nu=m^{\frac{n}{2}+1}\cdot \nu(V_{\bullet},m_{\bullet})=\nu_{v_{m}}(\Gamma)\; .
$$
Making the substitutions \eqref{substitutions} and using the relation $\gamma_{i}=\frac{r}{P}\Gamma_{i}$,
$$
\nu=m^{\frac{n}{2}+1}\cdot \frac {\sum_{i=1}^{t+1} \frac{\gamma_{i}}{r}[(r^i P - rP^i)]}
 {\sqrt{\sum_{i=1}^{t+1} P^{i}\frac{P^{2}}{r^{2}} \gamma_{i}^{2}}}\; ,
$$
whose square is a rational function on $m$ (since $P$ and $P^i$ are polynomials on $m$). Therefore we get
$$
\nu=m^{\frac{n}{2}+1}\cdot \frac{1}{P}\frac {\sum_{i=1}^{t+1}
\gamma_{i}[r^{i}P-rP^{i}]}
 {\sqrt{\sum_{i=1}^{t+1} P^{i}\gamma_{i}^{2}}}\; .
$$

\begin{prop}
\label{finalfunction} Given a sheaf $E$, there exists a unique filtration
$$0\subset E_{1} \subset E_{2} \subset \cdots \subset E_{t} \subset
E_{t+1}=E$$ with positive weights $n_{1},\ldots,n_{t}$, $n_{i}=\frac{\gamma_{i+1}-\gamma_{i}}{r}$, which gives
maximum for the function
$$\nu=m^{\frac{n}{2}+1}\cdot \frac {\sum_{i=1}^{t+1}
P^{i}\gamma_{i}[\frac{r^{i}}{P^{i}}-\frac{r}{P}]}
 {\sqrt{\sum_{i=1}^{t+1} P^{i}\gamma_{i}^{2}}}\; .$$
\end{prop}

Similarly, we had defined the coordinates $v_{i}$ (slopes of
segments of the graph), as
$$v_{i}=m^{n+1}\cdot \big[\frac{r^{i}}{P^{i}}-\frac{r}{P}\big]$$
Therefore we can express the function $\mu$ as
$$\nu=m^{-\frac{n}{2}}\cdot \frac {\sum_{i=1}^{t+1}
P^{i}\gamma_{i}v_{i}} {\sqrt{\sum_{i=1}^{t+1}
P^{i}\gamma_{i}^{2}}}=m^{-\frac{n}{2}}\cdot
\frac{(\gamma,v)}{||\gamma||}\; ,$$ where the scalar product is
given by the diagonal matrix
$$\left(
    \begin{array}{cccc}
      P^{1} &  &  & 0 \\
       & P^{2} &  &  \\
       &  & \ddots &  \\
      0 &  &  & P^{t+1} \\
    \end{array}
  \right)$$

\begin{prop}
\label{descendentslopes} Given the Kempf filtration of a sheaf
 $E$,
$$0\subset E_{1} \subset E_{2} \subset \cdots \subset E_{t} \subset E_{t+1}=E$$
it verifies
$$\frac{P^{1}}{r^{1}}>\frac{P^{2}}{r^{2}}>\ldots>\frac{P^{t+1}}{r^{t+1}}$$
\end{prop}
\begin{pr}
The coordinates of the vector $v$ associated to the filtration are, for $m$ large enough, $v_{i}=m^{n+1}\cdot
(\frac{r^{i}}{P^{i}}-\frac{r}{P})$. Now apply Lemma \ref{lemmaA} which says that $v$ is convex, i.e.
$v_{1}<\ldots<v_{t+1}$.
\end{pr}

\begin{prop}
\label{blocksemistability} Given the Kempf filtration of a sheaf $E$,
$$0\subset E_{1} \subset E_{2} \subset \cdots \subset E_{t} \subset E_{t+1}=E\; ,$$
each one of the blocks $E^{i}=E_{i}/E_{i-1}$ is semistable.
\end{prop}
\begin{pr}
Consider the graph associated to the Kempf filtration of $E$.
Suppose that any of the blocks has a destabilizing subsheaf. Then,
it corresponds to a point which lies above the graph of the filtration.
The graph obtained by adding this new point is a refinement of the
graph of the Kempf filtration, whose convex envelope is not the
original graph, which contradicts Lemma \ref{lemmaB}.
\end{pr}

\begin{cor}
The Kempf filtration of a sheaf $E$ coincides with its
Harder-Narasimhan filtration.
\end{cor}
\begin{pr}
By Propositions \ref{descendentslopes} and \ref{blocksemistability} the Kempf filtration verifies the two
properties of the Harder-Narasimhan filtration. By uniqueness of the Harder-Narasimhan filtration both
filtrations coincide.
\end{pr}


\begin{thebibliography}{EMG}

\bibitem[AB]{AB}{M. F. Atiyah and R. Bott}
\textit{The Yang Mills equations over Riemann surfaces, } Phil.
Trans. R. Soc. Lond. Ser. A \textbf{308} (1982), no. 1505,
523-615.

\bibitem[Br]{Br}{L. Bruasse, }
\textit{Optimal destabilizing vectors in some gauge theoretical moduli
  problems. }
Ann. Inst. Fourier (Grenoble) \textbf{56} (2006), no. 6, 1805-1826

\bibitem[BT]{BT}
L. Bruasse and A. Teleman
\textit{Harder-Narasimhan filtrations and optimal destabilizing
  vectors in complex geometry. }
Ann. Inst. Fourier (Grenoble) \textbf{55} (2005), no. 3, 1017-1053.

\bibitem[Gi]{Gi}{D. Gieseker, }
\textit{On the moduli of vector bundles on an algebraic surface, }
Ann. Math., \textbf{106} (1977), 45-60.

\bibitem[He]{He} {W. H. Hesselink}
\textit{Uniform instability in reductive groups, } J. Reine Angew.
Math. \textbf{304} (1978), 74-96.



\bibitem[HK]{HK}{V. Hoskins and F. Kirwan, }
\textit{Quotients of unstable subvarieties and moduli spaces of
sheaves of Fixed Harder-Narasimhan type, } Proc. Lond. Math Soc. (3) \textbf{105} (2012), 852-890.



\bibitem[HL0]{HL0}{D. Huybrechts and M. Lehn, }
\textit{Stable pairs on curves and surfaces, } J. Alg. Geom.,
\textbf{4} no. 1 (1995), 67-104.




\bibitem[HL1]{HL1}{D. Huybrechts and M. Lehn, } \textit{Framed
modules and their moduli, } Intern. J. Math., \textbf{6} no. 2
(1995), 297-324.

\bibitem[HL2]{HL2}{D. Huybrechts and M. Lehn, }
The geometry of moduli spaces of sheaves, Aspects of Mathematics
E31, Vieweg, Braunschweig/Wiesbaden 1997.


\bibitem[Ke]{Ke}{G. Kempf, }
\textit{Instability in invariant theory, } Ann. of Math. (2)
\textbf{108} no. 1 (1978), 299-316.


\bibitem[Ma]{Ma}{M. Maruyama, }
\textit{Moduli of stable sheaves, I and II.} J. Math. Kyoto Univ.
\textbf{17} (1977), 91--126. \textbf{18} (1978), 557-614.

\bibitem[GIT]{GIT}{D. Mumford, J Fogarty and F. Kirwan, }
Geometric invariant theory.
Third edition. Ergebnisse der Mathematik und ihrer
Grenzgebiete (2), 34. Springer-Verlag, Berlin, 1994.

\bibitem[Ne]{Ne}{P.E. Newstead, }
Lectures on Introduction to Moduli Problems and Orbit Spaces,
Published for the Tata Institute of Fundamental Research, Bombay.
Springer-Verlag, Berlin (1978).

\bibitem[RR]{RR}{S. Ramanan and A. Ramanathan, }
\textit{Some remarks on the instability flag, } T\^ohoku Math. Journ., \textbf{36} (1984),
269-291.

\bibitem[Si]{Si}{C. Simpson, }
\textit{Moduli of representations of the fundamental group of a
smooth projective variety I, } Publ. Math. I.H.E.S. \textbf{79}
(1994), 47--129.

\bibitem[Za1]{Za1}{A. Zamora, } \textit{Harder-Narasimhan
filtration for rank $2$ tensors and stable coverings, }
arXiv:1306.5651, (2013) (submitted preprint)

\bibitem[Za2]{Za2}{A. Zamora, }
\textit{On the stability of vector bundles, }
Master Thesis, Universidad Complutense de Madrid, 2009. Available at e-prints UCM server, 
http://www.mat.ucm.es/invesmat/wp-content/uploads/2011/12/trabajo-master-curso-2008-09-alfonso-zamora.pdf.


\bibitem[Za3]{Za3}{A. Zamora, }
\textit{GIT characterizations of Harder-Narasimhan filtrations, }
Ph.D. Thesis, Universidad Complutense de Madrid, 2013.

\end{thebibliography}
\end{document}